\documentclass[11pt,a4paper]{article}
\usepackage[utf8]{inputenc}
\usepackage[T1]{fontenc}
\usepackage{mathptmx}
\usepackage[dvipsnames]{xcolor}
\usepackage{microtype}
\usepackage[margin=2.5cm]{geometry}
\usepackage{amsmath,amssymb,amsthm,mathtools}
\usepackage{enumitem}
\usepackage[colorlinks=true,linkcolor=MidnightBlue,citecolor=MidnightBlue,urlcolor=MidnightBlue]{hyperref}

\theoremstyle{plain}
\newtheorem{theorem}{Theorem}
\newtheorem{proposition}{Proposition}
\theoremstyle{definition}
\newtheorem{definition}{Definition}

\newcommand{\PP}{\mathcal{P}}
\newcommand{\U}{\mathcal{U}}
\newcommand{\avg}[1]{\langle #1\rangle}
\newcommand{\avgU}[1]{\langle #1\rangle_{\U}}
\newcommand{\rhobar}{\bar\rho}
\newcommand{\ind}{\chi}
\newcommand{\Var}{\operatorname{Var}}
\newcommand{\nt}{\tilde n}

\title{\textbf{Invariant universal averages and Bertrand's paradox:\\ a response to Shackel}}
\author{Massimiliano Sassoli de Bianchi\\[6pt]
\small Center Leo Apostel for Interdisciplinary Studies, Brussels Free University, Brussels, Belgium}
\date{31 October 2026}

\begin{document}
\maketitle

\begin{abstract}
\noindent In 2014, with D.~Aerts we proposed a solution to the hard problem of Bertrand's paradox by means of a universal average over different ways of randomizing an entity. In his recent monograph, N.~Shackel raises objections concerning the approximation by cellular densities, the restrictions used in their construction, the apparently trivial limit of their averages, the choice of a second-order probability law, and the dependence of the result on the naming of chords. We revisit these objections using the invariance-based formulation of the universal average. The first objection is well taken, and the approximation argument is repaired: an explicit sequence of bivalent cellular densities provides weak approximation of any prescribed probability measure, with an error bound for intervals. We show that any second-order selection law invariant under rearrangements of equal-length cells has the same uniform barycenter, independently of the particular densities and weights employed. For cellular densities, this conclusion holds exactly at every finite resolution. Thus, non-uniqueness of the second-order law does not prevent uniqueness of the average it produces, which is what entitles one to speak of `the' universal average at all. Finally, we prove that no probability measure can be invariant under both these rearrangements and arbitrary changes of interval lengths. The reference geometry therefore cannot be obtained by averaging alone: it must be supplied by the physical modelization of the randomized entity. These results clarify the scope, and the limitations, of our original distinction between the easy and hard problems of Bertrand's paradox.
\end{abstract}

\section{Introduction}\label{sec:intro}

In \cite{AS2014b} we argued that Bertrand's question (``we draw at random a chord onto a circle; what is the probability that it is longer than the side of the inscribed equilateral triangle?'' \cite{Bertrand}) contains two distinct problems. The \emph{easy} problem is to specify a physical realization of the randomization, identifying the randomized entity and its geometric parameter space (for example, a stick, a pair of pebbles, or a pebble and a stick); the distinct specifications give rise to the classical alternative answers. The \emph{hard} problem is, once the entity is specified, to give a precise meaning to ``at random,'' i.e., to the condition of maximal lack of knowledge about \emph{the way} the entity is randomized. Our proposal for the hard problem was to describe the different ``ways of tossing'' the entity by probability densities $\rho$ on the relevant parameter space (for a stick, the lines meeting the circle, parameterized by their direction $\theta\in[0,\pi]$ and their distance $x\in[0,2]$ from the bottom of the circle), and to define a \emph{universal average} over all of them. Since the length of the chord depends only on $x$, the problem reduces to densities $\rho(x)$ on $[0,2]$, and the probability of a long chord in the $\rho$-way of tossing is 
\begin{equation}
\PP(L>\sqrt3|\rho)=\int_{1/2}^{3/2}\rho(x)\,dx
\end{equation}
The universal average was defined, in the spirit of the construction of the Wiener measure, by partitioning $[0,2]$ into $n$ cells of equal length, considering the $2^n-1$ \emph{cellular} densities taking a constant value on some of the cells and vanishing on the others, averaging the probabilities uniformly over them, and letting $n\to\infty$. We proved that this universal average coincides with the probability computed with the uniform density $\rho_u(x)=1/2$, which gives $1/2$ for the stick, and similarly $1/3$ for two pebbles and $1/4$ for a pebble inside the disk, using the uniform densities on the corresponding parameter spaces.

In his monograph on Bertrand's paradox \cite{Shackel2024}, N.~Shackel devotes an entire chapter (Ch.~9) to this construction, and returns to it in the following chapter on meta-indifference (Ch.~10). His critique is detailed and, in several respects, pertinent, and we are grateful to him for it. His objections (without claiming to be exhaustive) can be summarized as follows.

\begin{enumerate}[label=(O\arabic*),leftmargin=*]
\item \emph{Order of the quantifiers.} The argument by which we showed that every density $\rho$ can be approximated by cellular densities, in the sense that $\PP(x\in A|\rho_n)\to\PP(x\in A|\rho)$, as $n\to\infty$, only proves that for each interval $A$ there is a suitable sequence, not that one sequence works for all intervals, which is not sufficient for the claim that cellular densities are ``dense'' among all densities.
\item \emph{Three constraints.} The cellular densities are defined by three constraints, cells of equal length, two values only (one of which is zero), and equal weights. As soon as three or more values are allowed there are uncountably many cellular densities with $n$ cells, so that ``giving the same weight to each of them'' is no longer a finite average and seems undefinable.
\item \emph{Danger of triviality.} The average of the cellular densities of level $n$ is already the uniform density for every $n$, so that the universal average is the limit of a constant sequence, which suggests that the theorem is ``trivial.'' Moreover, the limits of sequences of bivalent cellular densities are said to be Dirac delta functions, so that the set $B_\omega$ of these limits would not be the set $K$ of all densities, and the construction would therefore not capture every possible probability law. Also, the interchange of the average and $n\to\infty$ limit is unjustified.
\item \emph{Vicious regress.} Any average over the set of probability measures requires a second-order measure (Shackel's Principle of Meta-indifference for Sets, and his Theorem~40). There are infinitely many ways of constructing such measures, through the various norms on function spaces or otherwise, there is no measure that makes sense as \emph{the} uniform measure, and different choices give contradictory values, so that Bertrand's paradox recurs at the second order, and a retreat to a third-order principle leads to a vicious regress.
\item \emph{Namings.} Our parameterization of the stick is a ``naming'' of the chords, and the universal average can be defined for any naming; but different namings give different values, hence the universal average does not avoid the recurrence of the paradox, and our approach is, after all, an instance of the Distinction strategy, which Shackel regards as falling into van Fraassen's trap.
\end{enumerate}

\noindent The purpose of this note is to reply to these objections. The reply draws on a companion article \cite{Sassoli2026}, in which the notion of universal average is re-examined in its original context, that of the hidden-measurement interpretation of quantum mechanics \cite{AS2014a}, and in which the universal average is shown to be an \emph{invariant average}. More precisely, in Section~\ref{sec:approx}, we repair the approximation argument (O1), giving an explicit single sequence with an error estimate valid for all intervals, and weak convergence to any prescribed probability law. In Section~\ref{sec:invariant}, we show that (O2) and (O3) are dissolved once the universal average is understood as an invariant average: its value is fixed by symmetry, exactly, at every finite resolution, for every selection of densities which does not privilege any portion of the parameter space, whatever the number of values, the weights, or the limit. In Section~\ref{sec:meta}, we address (O4): Shackel's meta-indifferent average has the same mathematical form as our universal average, and the theorem of Section~\ref{sec:invariant} shows that all second-order laws satisfying our explicit rearrangement-indifference condition yield the same barycenter, without requiring a canonical choice among those laws. In other words, there is no vicious regress. In Section~\ref{sec:naming}, we address (O5), and we concede its mathematical core, in the form of a theorem: no notion of size can be extracted from the set of all densities. We argue that this does not refute but confirms the distinction between the easy and the hard problem, which it shows to be unavoidable.

\section{The approximation argument}\label{sec:approx}

Throughout, we retain the term ``density'' in the broad sense used in \cite{AS2014b}: when Dirac delta or singular components are present, $\rho$ denotes the corresponding probability measure, and expressions such as $\int_A\rho(x)\,dx$ mean the probability assigned to $A$.

We first recall the construction of \cite{AS2014b}. The interval $[0,2]$ is partitioned into $n=m\ell$ elementary cells of length $2/n$, grouped into $m$ cells $S_1,\dots,S_m$ of length $2/m$, each containing $\ell$ elementary cells. A (bivalent) cellular density is
\begin{equation}\label{eq:cellular}
\rho_n(x)=\frac{n}{2\nt}\,\chi_{m\ell}(x)
\end{equation}
where $\chi_{m\ell}$ is a step function taking the value $1$ in $\nt\ge1$ elementary cells (the non-zero cells) and $0$ in the others. The claim of Section~VI of \cite{AS2014b} was that for every density $\rho$ on $[0,2]$ a single sequence of cellular densities $\rho_n$ can approximate its interval probabilities. For targets with point masses, the natural convergence statement must be restricted to intervals whose boundaries carry no target mass. Shackel is right that the argument originally given did not establish even this simultaneous convergence. It chose the number $\nt_i$ of non-zero elementary cells in each cell $S_i$ so that $\nt_i/\nt\to\int_{S_i}\rho(x) dx$, using the density of the rationals in $[0,1]$, cell by cell, so that it only established that for each interval there is a suitable sequence.
Here is the correct statement, with its proof.

\begin{proposition}\label{prop:approx}
Let $\rho$ be an arbitrary probability measure on $[0,2]$, expressed in the generalized-density notation just introduced. There is a single sequence of cellular densities $\rho_n$, $n=m\ell$, $\ell=m^2+m$, $m=1,2,\dots$, such that
\begin{equation}\label{eq:bound}
\big|\PP(x\in A|\rho_n)-\PP(x\in A|\rho)\big|<\frac4m+4\max_{1\le i\le m}\int_{S_i}\rho(x) dx
\end{equation}
for every interval $A\subset[0,2]$, where $S_i$ are the $m$ cells of level $m$. Consequently, as $m\to\infty$ (and hence $n\to\infty$), the corresponding probability measures converge weakly to $\rho$: their interval probabilities converge whenever the endpoints carry zero $\rho$-mass, and they converge uniformly over all intervals if $\rho$ has no point masses.
\end{proposition}

\begin{proof}
Let $p_i=\int_{S_i}\rho(x)dx$, $i=1,\dots,m$, be the cell masses of $\rho$, $\sum_ip_i=1$, and let $\rho_n$ have, in the cell $S_i$, exactly $\nt_i=\lfloor\ell p_i\rfloor$ non-zero elementary cells (the first $\nt_i$ ones, say), so that the total number of non-zero elementary cells is $\nt=\sum_i\nt_i$. Since each $\lfloor\ell p_i\rfloor$ differs from $\ell p_i$ by less than $1$, we have $\ell-m<\nt\le\ell$; in particular $\nt\ge1$, so that $\rho_n$ is well defined. The mass of $\rho_n$ in $S_i$ is $\nt_i/\nt$, and the total discrepancy between the cell masses of $\rho_n$ and those of $\rho$ is
\begin{equation}\label{eq:tv}
\sum_{i=1}^m\Big|\frac{\nt_i}{\nt}-p_i\Big|\le\sum_{i=1}^m\frac{|\nt_i-\ell p_i|}{\nt}+\sum_{i=1}^m\ell p_i\Big(\frac1{\nt}-\frac1\ell\Big)<\frac{m}{\nt}+\frac{\ell-\nt}{\nt}<\frac{2m}{\ell-m}=\frac2m
\end{equation}
where we used $\ell-\nt<m$ and $\nt>\ell-m$ in the last inequality. Now let $A$ be an interval whose endpoints fall in the cells $S_a$ and $S_b$ (possibly $a=b$), and let $U$ be the union of the cells entirely contained in $A$ ($U$ may be empty). By \eqref{eq:tv}, $|\PP(x\in U|\rho_n)-\PP(x\in U|\rho)|<\frac2m$. The part of $A$ outside $U$ is contained in $S_a\cup S_b$, which carries a mass at most $p_a+p_b$ under $\rho$, and at most $p_a+p_b+\frac2m$ under $\rho_n$, again by \eqref{eq:tv}. Altogether, we obtain
\begin{equation}
\big|\PP(x\in A|\rho_n)-\PP(x\in A|\rho)\big|<\frac2m+\Big(p_a+p_b+\frac2m\Big)+(p_a+p_b)\le\frac4m+4\max_ip_i
\end{equation}
which is \eqref{eq:bound}. If $\rho$ has no delta components, $\max_ip_i\to0$ as $m\to\infty$ (a measure without atoms gives vanishing mass to cells of vanishing length), and the right-hand side tends to zero uniformly in $A$. For a fixed interval $A$ whose endpoints have zero $\rho$-mass, the masses $p_a$ and $p_b$ of the shrinking cells containing those endpoints tend to zero even if $\rho$ has atoms elsewhere. The preceding interval estimate therefore also proves convergence on every such interval, hence weak convergence.
\end{proof}

\noindent Two additional remarks are in order. First, the uniqueness of a measure determined by its values on intervals must be distinguished from the convergence of a sequence of measures. Intervals form a $\pi$-system generating the Borel $\sigma$-algebra, and two probability measures agreeing on such a system agree on all Borel sets (Dynkin's $\pi$--$\lambda$ theorem; see, e.g., \cite[Theorem~3.3]{Billingsley1995}). This uniqueness result does \emph{not} imply that probabilities converge on every Borel event: for example, cellular densities cannot converge to $\delta_z$ on the singleton $\{z\}$, which always has zero probability under a genuine density. The appropriate statement is weak convergence, equivalently convergence on all continuity sets of the limiting measure, by the portmanteau theorem \cite[Theorem~2.1]{Billingsley1999}. Proposition~\ref{prop:approx} establishes this with an explicit construction and an error bound for intervals. Shackel's observation concerning the distinction between Borel and Lebesgue $\sigma$-algebras does not affect this claim, which is formulated for Borel probability measures. No identification of arbitrary completions with Lebesgue measurable sets is required.

Second, Shackel's Theorem~39 spells out what the density of the set of all cellular densities would mean in the \emph{supremum} metric, i.e., uniform convergence of the functions $\rho_n$ to $\rho$, and he observes that we did not prove it. Indeed, it is false for bivalent (two-valued) step functions, but nobody needs it: what is physically meaningful is not the pointwise value of a density but the probabilities it assigns to events, and for these Proposition~\ref{prop:approx} gives what is needed. (For an absolutely continuous target density $\rho$, dropping the two-value restriction and using its cell-average histograms yields $L^1$ convergence, hence convergence uniformly over all Borel events, by a standard result of analysis; see Section~\ref{sec:invariant}. This stronger claim need not hold for targets with point masses or singular-continuous parts). As we shall now see, however, in the correct formulation of the universal average the approximation property is no longer what the computation of the average rests upon.

\section{The universal average as an invariant average}\label{sec:invariant}

The construction of \cite{AS2014b} averaged over a specific finite family of densities at each level, and then took a limit. The objections (O2) and (O3) are both about the status of this family and of this limit. The answer to both is that neither the family nor the limit is essential. What is essential is a \emph{symmetry}, and once it is identified the universal average is fixed by it, exactly and uniquely.

\subsection{Indifferent selections}

A ``way of tossing'' is a density $\rho$. An average over the ways of tossing is defined by a probability law $\U$ on the set of densities, describing a random selection of $\rho$. 
The average of the probability of an event $A$ is then
\begin{equation}\label{eq:avgU}
\avgU{\PP(x\in A|\rho)}=\int\PP(x\in A|\rho)\,d\U(\rho)
\end{equation}
The construction of \cite{AS2014b} corresponds, at each level $n$, to the law $\U_n$ giving probability $1/(2^n-1)$ to each of the $2^n-1$ cellular densities \eqref{eq:cellular}. Since $\PP(x\in A|\rho)$ depends linearly on $\rho$, the average of the probabilities is computed with the \emph{average probability measure} (or barycenter) $\rhobar_\U\equiv\int\rho\,d\U(\rho)$, whose density will be uniform in the indifferent case:
\begin{equation}\label{eq:key}
\avgU{\PP(x\in A|\rho)}=\PP(x\in A|\rhobar_\U)=\int_A\rhobar_\U(x)\,dx
\end{equation}
The whole question is therefore: for which laws $\U$ is $\rhobar_\U$ the uniform probability measure, with density $\rho_u=1/2$?

Consider the $n$ cells $\sigma_j^{(n)}=[(j-1)\frac2n,\,j\frac2n)$, $j=1,\dots,n$, of equal length (the last one closed at $x=2$), and the cell masses $m_j(\rho)=\int_{\sigma_j^{(n)}}\rho(x)dx$, $\sum_jm_j(\rho)=1$. A \emph{rearrangement} of the interval is the operation of cutting it along the cell boundaries and gluing the $n$ pieces again in a different order, specified by a permutation $\pi$: it maps a density $\rho$ to a density $\rho^\pi$ whose cell masses are permuted, $m_{\pi(j)}(\rho^\pi)=m_j(\rho)$, and whose profile inside each cell is the profile that $\rho$ had in the piece moved there. A rearrangement preserves lengths: it maps genuine densities to genuine densities, delta components to delta components, and the uniform density to itself.

\begin{definition}[indifferent selection]\label{def:indiff}
A law $\U$ on the set of densities is \emph{indifferent} if, for every $n$ of the form $2^k$ and every permutation $\pi$ of the $n$ cells, the law of the rearranged random density $\rho^\pi$ is the same as the law of $\rho$, i.e., $\U(\rho^\pi\in E)=\U(\rho\in E)$ for every (measurable) set $E$ of densities. In particular, the joint law of the cell masses $(m_1(\rho),\dots,m_n(\rho))$ is symmetric under permutations of the cells. 
The average \eqref{eq:avgU} produced by an indifferent law is what we shall call a \emph{universal average}.
\end{definition}

\noindent Note that the restriction to $n=2^k$ is just a convenience: dyadic cells are nested, but any nested sequence of partitions into cells of equal length would also do. Indifference is the precise form of the requirement that the selection of a way of tossing ``does not privilege any portion of the parameter space'': selecting $\rho$ or selecting its rearranged version $\rho^\pi$ is, for an indifferent law, the same thing. It is the second-order counterpart of the indifference that, at the first order, selects the uniform density: a condition of maximal lack of knowledge about \emph{the way} the chord is produced, rather than about the chord itself. 
This is why we reserve the term \emph{universal} for the averages that indifferent laws produce, and not for every average over ways of tossing: a law concentrated on a single density averages nothing, and returns the probability given by that density alone.\footnote{This sharpens the usage of \cite{AS2014b,AS2014a}, where the universal average was identified with the particular limiting construction recalled above; see \cite{Sassoli2026} for a fuller discussion of the terminology.} The law $\U_n$ of \cite{AS2014b} is indifferent at its own level $n$ (a rearrangement maps a cellular density with $\nt$ non-zero cells to another one with the same $\nt$), though not at the finer levels $2n,4n,\dots$, since the two halves of a non-zero cell always carry the same mass; it therefore does not satisfy Definition~\ref{def:indiff} as it stands, and Theorem~\ref{thm:finite} below will show that it nevertheless yields the universal average exactly. But the definition refers to no family, no discretization and no limit: it is a property of $\U$ alone.

\subsection{Uniqueness}

\begin{theorem}[uniqueness of the universal average]\label{thm:unique}
If $\U$ is indifferent, then its barycenter is $\rhobar_\U=\rho_u$, the uniform probability measure on $[0,2]$, with density $1/2$. Consequently, for every event $A$,
\begin{equation}
\avgU{\PP(x\in A|\rho)}=\PP(x\in A|\rho_u)=\frac{|A|}{2},
\end{equation}
and in particular $\avgU{\PP(L>\sqrt3|\rho)}=\frac12$, whatever the indifferent law $\U$.
\end{theorem}

\begin{proof}
Fix a dyadic level $n$. The cell masses of the barycenter are the averages of the cell masses of $\rho$: $m_j(\rhobar_\U)=\avgU{m_j(\rho)}$. By indifference, for any two cells $j$ and $l$, the random variables $m_j(\rho)$ and $m_l(\rho)$ take their values with the same probabilities [exchange the two cells by a rearrangement $\pi$: then $m_l(\rho^\pi)=m_j(\rho)$, and $\rho^\pi$ has the same law as $\rho$], so that they have the same mean value. The $n$ numbers $m_j(\rhobar_\U)$ are thus all equal, and since they sum to $1$, each equals $1/n$, which is the mass that $\rho_u$ gives to a cell of length $2/n$. Hence, the cumulative function $\bar F(x)=\int_0^x\rhobar_\U(y)dy$ coincides with $F_u(x)= x/2$ at every dyadic point $x=2j/2^k$. Both functions being non-decreasing and $F_u$ continuous, squeezing an arbitrary $x$ between consecutive dyadic points of level $k$ and letting $k\to\infty$ gives $\bar F=F_u$ everywhere, i.e., $\rhobar_\U=\rho_u$ (as probability measures, hence as densities almost everywhere).
\end{proof}

\noindent In other words, if the average density were not uniform, some cell would receive an average mass different from that of another cell of the same length, whereas by hypothesis the selection does not distinguish between the two. The proof uses only the symmetry of $\U$ and the linearity of the probabilities in $\rho$. It is insensitive to the kind of densities selected (two-valued, many-valued, continuous, with or without delta components), to their number, to the weights (which must be equal only for densities that are rearrangements of one another), to the resolution, and to the existence of a limit. The universal average is thus \emph{unique}: it does not depend on the way the ``average over all $\rho$'' is implemented, as long as the implementation is indifferent, and we may define, independently of any construction,
\begin{equation}\label{eq:univ}
\avg{\PP(x\in A|\rho)}^{\rm univ}\equiv\avgU{\PP(x\in A|\rho)},\quad\U\text{ any indifferent law}
\end{equation}
It is Theorem~\ref{thm:unique}, and not a convention, that entitles us to the definite article in ``\emph{the} universal average'': the indifferent laws are infinitely many, but they all assign the same probability to every event, so that for the purpose of computing probabilities they are interchangeable. This is the precise sense in which the construction of \cite{AS2014b} was one realization of a unique object, rather than one choice among incompatible alternatives -- which is what the objections (O2) and (O4) suspect it to be.

\subsection{Finite resolution}

The bivalent cellular laws $\U_n$ of \cite{AS2014b} are indifferent at a single level. More generally, call \emph{cellular densities of level $n$} the densities that are constant in each of the $n$ cells, i.e., $\rho_p=\frac n2\sum_jp_j\ind_{\sigma_j^{(n)}}$, with arbitrary cell masses $p=(p_1,\dots,p_n)$, $p_j\ge0$, $\sum_jp_j=1$, and call $\tilde\U_n$ any law on them. The bivalent densities \eqref{eq:cellular} are those with $p_j\in\{0,1/\nt\}$. A rearrangement maps $\rho_p$ to $\rho_{p^\pi}$, with $p^\pi_{\pi(j)}=p_j$, and the orbit of $\rho_p$ consists of the densities $\rho_{p^\pi}$ obtained through the $n!$ permutations (some of which may coincide). We have the following theorem:

\begin{theorem}[finite resolution]\label{thm:finite}
The average of every orbit is the uniform density, $\frac{1}{n!}\sum_\pi\rho_{p^\pi}=\rho_u$. Consequently, every law $\tilde\U_n$ on the cellular densities of level $n$ which is invariant under rearrangements has $\rhobar_{\tilde\U_n}=\rho_u$ exactly, with no limit to be taken. In particular, this holds for the bivalent laws $\U_n$, for every $n$ and every event $A$.
\end{theorem}

\begin{proof}
A rearrangement applied to the orbit average merely permutes the $n!$ terms of its sum, so the orbit average is a cellular density invariant under all rearrangements. Its cell masses are therefore all equal, hence equal to $1/n$, so that it coincides with $\rho_u$. 
If $\tilde\U_n$ is invariant, $\avg{\rho}_{\tilde\U_n}=\avg{\rho^\pi}_{\tilde\U_n}$ for every $\pi$. Averaging over $\pi$ gives $\avg{\rho}_{\tilde\U_n}=\avg{\rho_u}_{\tilde\U_n}=\rho_u$.
\end{proof}

\subsection{Replies to (O2) and (O3)}

Of the three constraints (equal length, two values and equal weights), only the equal length of the cells matters, and it matters because it encodes indifference. The two values and the equal weights are instead conveniences. Shackel is right that with three or more values there are uncountably many cellular densities of level $n$, but there is no need to ``give the same weight to each of them'' as if they were finitely many. One averages them with any law invariant under rearrangements (for instance the uniform measure on the simplex of the cell-mass vectors, or any symmetric Dirichlet law), and Theorem~\ref{thm:finite} gives exactly $\rho_u$. The weights need not be uniform either: any weighting of the orbits gives $\rho_u$. The combinatorial identities (24) and (29) of \cite{AS2014b}, used in the proof, are simple consequences of the symmetry; see \cite{Sassoli2026}, Appendix~A.

Coming to the sequence of averages, it is constant because the symmetry already fixes the average at every finite resolution. This is not a defect of the result but its content: the value is exact, independent of the resolution and of the family, as the center of mass of a symmetric body lies on its axis of symmetry independently of how finely the body is subdivided. The limit $n\to\infty$ in \cite{AS2014b} added nothing to the \emph{value} of the average. What it was meant to secure is that no density is left out, i.e., that the family of levels is exhaustive. This is a different matter, to which we now turn.

Shackel challenges the identification of the set $B_\omega$ of limits of the chosen two-valued sequences with the set $K$ of all densities, as well as the unjustified interchange between an average and a limit. The criticism of the original limiting argument is pertinent, but it does not imply that every pointwise limit of bivalent densities is a Dirac mass, as he mentions (the constant uniform sequence is an immediate counterexample). Neither point is needed in the present formulation. What is relevant is not the set of pointwise limits of sequences of cellular densities, but the fact that the \emph{union} of the levels is weakly dense in $K$ (in the sense of Proposition~\ref{prop:approx} for two-valued densities, and in the $L^1$ sense for histograms). And no interchange of limit and average occurs anywhere: by Theorems~\ref{thm:unique} and \ref{thm:finite}, the average is $\rho_u$ at each level and for every indifferent law, so that there is simply nothing to interchange. The identification of the universal average with an indifferent average over ways of tossing does not rest on a limiting procedure but on the uniqueness theorem: every indifferent way of averaging over the densities, finite or infinite, gives the same answer, and the next section gives an indifferent law whose weak support contains every probability measure on the interval.

\section{Meta-indifference}\label{sec:meta}

In Ch.~10, Shackel formalizes meta-indifference by joint averaging as follows [his Principle of Meta-indifference for Sets (POMS) and Theorem~40]. Using his notation, let $\Psi$ be the set of first-order probability measures on the event space. Equip $\Psi$ with a $\sigma$-algebra and a finite measure $\nu$ expressing ``equal epistemic status'' (members of equal $\nu$-measure have equal second-order probability). Then the meta-indifferent probability measure is the expectation of $P(\sigma)$ over the normalized second-order measure:
\begin{equation}\label{eq:PPsi}
P_\Psi(\sigma)=\int_\Psi V_\sigma\,d\Big(\frac{\nu}{\nu(\Psi)}\Big),\quad V_\sigma(P)=P(\sigma)
\end{equation}
This is exactly our universal average: with $\nu/\nu(\Psi)=\U$, \eqref{eq:PPsi} is \eqref{eq:avgU}, and $P_\Psi$ is the measure with density $\rhobar_\U$. We are therefore in complete agreement with Shackel on \emph{what} a meta-indifferent average is. The disagreement is on whether it is well defined. Shackel's argument is that there are infinitely many ways of putting a measure on $\Psi$ (through the metrics induced by the $L^n$ norms, through pre-measures, or otherwise), and since (p.~234) ``there is no measure which makes sense as \emph{the} uniform measure,'' different choices give contradictory values of $P_\Psi$. Hence the paradox recurs at the second order, and a third-order principle would start a vicious regress.

Theorem~\ref{thm:unique} answers this argument directly. It is true that there is no canonical uniform measure on $\Psi$, and we shall see in a moment that the situation is even worse than Shackel says. But a canonical measure is not needed. POMS calls for equal second-order weight to be assigned to possibilities judged epistemically equivalent. In the present model we make this equivalence precise by requiring invariance under rearrangements of equal-length cells: these rearrangements preserve the reference geometry, so the selection of a way of tossing does not privilege one portion of the parameter space over another. This is not a further arbitrary choice added on top of POMS, but the natural operational reading of ``equal epistemic status'' once the reference geometry is fixed: two ways of tossing that differ only by which piece of the parameter space sits where are distinguished by nothing available in a condition of maximal lack of knowledge, and a law treating them differently would have to appeal to information that, by hypothesis, is absent. With this operational interpretation of meta-indifference, the normalized second-order law $\nu/\nu(\Psi)$ is indifferent in the sense of Definition~\ref{def:indiff}, and Theorem~\ref{thm:unique} says that \emph{all} such laws give the same $P_\Psi$, namely the uniform measure. The choice among them, which indeed cannot be made canonically, does not need to be made. There is no recurrence of the paradox at the second order, and no third-order principle is called for: the regress stops at the first step, because all the admissible second-order measures agree on the only thing that matters, their barycenter. A second-order measure constructed from an $L^n$ metric, or in any other way, either satisfies the stated rearrangement invariance, and then it gives the uniform answer, or it does not, and then it assigns different second-order weights to ways of tossing that differ only by the location of their features along the parameter space. The general principle of meta-indifference does not by itself exclude the second possibility; but a law of that kind does privilege a portion of the parameter space over another of the same size, which is what meta-indifference was introduced to avoid. Shackel's example involving the three processes of question BQ4 \cite{AS2014b}, with equal weights giving $\frac13(\frac12+\frac13+\frac14)=\frac{13}{36}$, illustrates instead a distinct randomization that first selects among different types of entities. The resulting answer is consistent with that specified protocol; it does not supply a competing indifferent law for a single fixed geometrical model.

Two further results of \cite{Sassoli2026} complete the picture, and we state them here without proof. First, Shackel is right that there is no uniform measure on $\Psi$, in the following strong sense: \emph{there is no non-trivial indifferent law supported on measures without point masses}: the only such law is concentrated on $\rho_u$ itself. More precisely, for an indifferent $\U$, the variance of every cell mass is $\Var_\U[m_j(\rho)]=\beta\frac1n(1-\frac1n)$, where $\beta=\avgU{\sum_iw_i^2}$ is the average of the squared weights of the delta components of the selected $\rho$. So, if $\U$ selects measures without point masses only, $\beta=0$, all variances vanish, and $\rho=\rho_u$ with probability one. This identifies a limitation of attempts to construct a non-trivial indifferent selection law while restricting attention to non-atomic measures: within that class only the trivial uniform selection exists. The finite-family construction of \cite{AS2014b} bypassed the need to specify such a law. Second, non-trivial indifferent laws exist on the set of \emph{all} probability measures (including those with delta components), the \emph{Dirichlet processes} of probability theory \cite{Ferguson} with base measure $\rho_u$ being the standard example. They are indifferent, and, with the uniform base measure, their \emph{weak support} is the set of all Borel probability measures on $[0,2]$ \cite{Ferguson,Majumdar1992}: every weak neighborhood of any prescribed measure has strictly positive probability. Equivalently, the probability masses of any prescribed finite partition into intervals of positive length can be approximated with positive probability, even for targets containing point masses. By Theorem~\ref{thm:unique} their barycenter is $\rho_u$. Hence a single indifferent second-order law has all probability measures in its topological support, without selecting each one exactly, and its average is again uniform.

\section{Easy and hard problems}\label{sec:naming}

We come to (O5), which is, we agree with Shackel, the crux of the matter, and on which our position has become more precise since \cite{AS2014b}.

\subsection{What is correct in the objection}

Shackel observes that our parameterizations are namings of the chords, that the universal average can be defined for any naming, and that different namings give different values. The mathematical core of this observation is correct, and it can be stated as a theorem. Consider, besides the length-preserving rearrangements of Section~\ref{sec:invariant}, the following \emph{generalized rearrangements}: cut the parameter space into pieces of \emph{un}equal lengths, and glue them again in a possibly different order after stretching or compressing each piece to fit its new place. These operations, understood as bijective piecewise-affine maps of the interval together with the push-forward of probability measures (and the corresponding Jacobian for ordinary densities), map the set $K$ of probability measures onto itself, just as the length-preserving rearrangements do. A coordinate change $y=h(x)$ alone preserves all probabilities when the reference measure is transformed consistently. Instead, declaring $dy$ to be a new uniform reference measure changes the randomization structure.

\begin{theorem}[no size from the set of densities]\label{thm:size}
No Borel probability measure on $[0,2]$ is invariant under all generalized rearrangements, including the length-preserving ones. Consequently, no second-order law $\U$ can be indifferent with respect to all of them: the universal average cannot be defined from the set $K$ of all densities alone, without a notion of size on the parameter space.
\end{theorem}

\begin{proof}
Assume that a probability measure is invariant under all the stated transformations. It is then invariant, in particular, under all length-preserving dyadic rearrangements. By the cell-mass argument of Theorem~\ref{thm:unique}, its mass on every dyadic cell of level $n$ equals $1/n$, so the probability measure is the uniform one on $[0,2]$. Consider now the piecewise-affine, length-changing bijection $h:[0,2]\to[0,2]$ defined by $h(x)=x/2$ on $[0,1]$ and $h(x)=3x/2-1$ on $[1,2]$. Its push-forward sends the mass $1/2$ of $[0,1]$ into $[0,1/2]$, whereas the uniform measure assigns mass $1/4$ to $[0,1/2]$. Thus the uniform measure is not invariant under $h$, a contradiction. If a second-order law $\U$ were invariant under all these transformations, its barycenter would also be invariant, giving the same contradiction.
\end{proof}

\noindent So, indifference, at any order, is always indifference \emph{with respect to a notion of size}, and no averaging procedure can produce this notion out of the bare set of densities. Shackel's ``Deluge theorem'' and his analysis of namings express, in his framework, the same fact. We fully concede it.

\subsection{What the objection does not show}

What Theorem~\ref{thm:size} does \emph{not} say is that the notion of size is arbitrary. In our approach it is supplied by the solution of the easy problem, and this is not a choice made at the second order, after the fact, but a datum of the modelization. The stick is a straight line thrown onto the plane, and the parameter space $[0,2]\times[0,\pi]$ with the measure $dx\,d\theta$ is not just ``a naming'': $dx\,d\theta$ is the \emph{invariant measure} on the lines of the plane, the unique (up to a constant factor) measure on sets of lines that is invariant under the group of rigid motions (translations and rotations) of the plane (\cite[Ch.~3, Sec.~1, pp.~27--29]{Santalo}), i.e., the measure determined by the symmetry of the randomized entity.\footnote{\emph{Kinematic measure} is Santal\'o's name for the invariant measure on the group of motions itself (Ch.~6). For lines he speaks of the invariant density $dG=dp\wedge d\theta$, with $p$ the signed distance from the origin (when $\theta\in[0,\pi)$), which is our $dx\,d\theta$ up to the shift $x=p+1$.} 
This is also why our value for the stick, $1/2$, coincides with the value obtained by Jaynes \cite{Jaynes} by requiring invariance under rotations, translations and dilations: Jaynes imposes the symmetry on the answer, we impose it, through indifference, on the ways of tossing, and Theorem~\ref{thm:unique} says that the two agree. 

For two pebbles tossed independently onto the circumference, the parameter space is the torus $[0,2\pi]^2$ with its product measure, the Haar measure of the independent rotations of the two points. For a pebble inside the disk, the parameter space is the disk itself, with its area measure. 
Shackel's quotient naming of the pebble case (unordered pairs of endpoints) is the push-forward of this second modelization, and it gives $1/3$. The $(x,\theta)$ naming carries the symmetry of the lines and gives $1/2$. The two namings name the same chords, but they describe two \emph{different randomized entities}, which was the point of \cite{AS2014b}. A mere change of naming does not alter the randomized entity or its probabilities, provided the underlying geometric measure is transformed consistently. What changes the problem is adopting a different reference measure, for example by declaring the new coordinates uniformly distributed. Such a choice amounts to specifying a different randomization structure, and therefore belongs to the easy problem rather than to the calculation of an average at fixed geometry.

This leads us to a conclusion which is, we believe, stronger than the one we drew in \cite{AS2014b}, and which answers Shackel's claim that our approach is ``merely'' an instance of the Distinction strategy. Bertrand's question asks for a chord chosen ``at random,'' i.e., it appeals to indifference, and Theorem~\ref{thm:size} shows that indifference, however it is implemented (by the principle of insufficient reason applied to the chords, by a universal average over the ways of selecting them, or by any other procedure that does not privilege some chords over others) cannot produce a notion of size out of nothing. Consequently, \emph{no} solution based on indifference can avoid the separation between an easy problem, which specifies the physical realization of the randomization and thereby its reference geometry, and a hard problem, which, given that specification and the adopted indifference condition, determines the resulting probabilities. This distinction is therefore a necessary part of an indifference-based solution. Shackel's observation that POMS and his Theorem~40 ``cannot fail to entail a Meta-indifferent Probability Measure'' (Ch.~10, p.~234) is compatible with this result: once the geometry and the explicit rearrangement-indifference condition are specified, joint averaging entails a unique resulting probability measure.

His further argument that this probability is not unique because the second-order measure is not unique fails, by the same theorem. What remains non-unique is the geometry, and this non-uniqueness lives in the easy problem, where Bertrand's question, as posed, leaves it. Whether one prefers to say that the paradox is thereby solved, or that Bertrand's question is shown to be incomplete, is largely a matter of taste. What is not a matter of taste is that the hard problem, once the easy one is solved, has a unique, exact and robust solution.

As for van Fraassen's trap, i.e., the worry that distinguishing questions by their physical realizations leads to as many answers as there are realizations, we can only observe that this is a feature of the question, not of the method: an ill-posed question has as many answers as there are ways of completing it. What the universal average contributes is to show that, once a physically motivated reference geometry and its rearrangement-indifference condition have been specified, no further choice of a particular second-order law is needed to fix the resulting average.

\section{Conclusion}\label{sec:conclusion}

We have replied to Shackel's objections to the universal average of \cite{AS2014b}. The technical objection concerning the order of the quantifiers was pertinent, and Proposition~\ref{prop:approx} now establishes the appropriate weak-approximation result through a single explicit sequence. The objections concerning the arbitrariness of the cellular densities, the constancy of the sequence of averages, and the interchange of limit and average are dissolved by Theorems~\ref{thm:unique} and \ref{thm:finite}: within the explicitly defined rearrangement-indifferent class, the barycenter is fixed by symmetry and unique, which is what entitles us to speak of \emph{the} universal average; for cellular densities the result is exact at every finite resolution. 
This also answers the objection concerning meta-indifference by joint averaging: Shackel's $P_\Psi$ has the form of our universal average, while the non-uniqueness of the second-order law is immaterial to the barycenter within the specified indifferent class. Finally, the objection concerning namings is correct in its mathematical core, which we have stated as Theorem~\ref{thm:size}: no notion of size can be extracted from the set of densities. But this confirms, rather than refutes, the distinction between the easy and the hard problem of Bertrand's paradox. The notion of size is supplied by the modelization of the randomized entity, through its symmetry, and once it is supplied the hard problem has a unique solution. The paradox arises, and will keep arising, whenever the easy problem is left implicit.

\section*{Acknowledgments}
I am grateful to Nicholas Shackel for the care with which he examined our work, and for the correspondence that preceded his book.

\section*{Declaration of generative AI use}
This reply draws on results presented in \cite{Sassoli2026}. Claude (Anthropic) assisted with the organization and drafting of the original manuscript and with scientific checking. ChatGPT (OpenAI) was subsequently used to assist with mathematical and bibliographical checking, critical review, and revision of the exposition. These tools were employed to identify potential weaknesses, test formulations, and improve clarity; they did not replace the author's responsibility for independently checking the arguments and sources. The author takes full responsibility for the final content.

\end{document}